\documentclass[12pt]{article}

\usepackage{amsmath}
\usepackage{amssymb}
\usepackage[mathscr]{eucal}
\usepackage{graphicx}
\usepackage{color}

\newtheorem{Theorem}{\sc Theorem}
\newtheorem{Definition}[Theorem]{\sc Definition}
\newtheorem{Proposition}[Theorem]{\sc Proposition}

\newtheorem{Corollary}[Theorem]{\sc Corollary}

\newtheorem{Example}[Theorem]{\sc Example}

\newcommand{\cL}{\mbox{{${\cal L}$}}}
\newcommand{\cQ}{\mbox{{${\cal Q}$}}}
\newcommand{\cM}{\mbox{{${\cal M}$}}}
\newcommand{\cP}{\mbox{{${\cal P}$}}}

\newcommand{\cS}{\mbox{{${\cal S}$}}}

\newcommand{\R}{{\if mm {\rm I}\mkern -3mu{\rm R}\else \leavevmode
\hbox{I}\kern -.17em\hbox{R} \fi}}

\newcommand{\cT}{\mbox{{${\cal T}$}}}
\newcommand{\cC}{\mbox{{${\cal C}$}}}

\newcommand{\bx}{\mbox{\boldmath{$x$}}}

\newcommand{\bnu}{\mbox{\boldmath{$\nu$}}}

\newcommand{\cH}{\mbox{{${\cal H}$}}}
\newcommand{\wu}{\mbox{{$\widetilde{u}$}}}

\newcommand{\wK}{\mbox{{$\widetilde{K}$}}}

\def\sqr#1#2{{
    \vcenter{
         \vbox{\hrule height.#2pt
               \hbox{\vrule width.#2pt height#1pt \kern#1pt
                     \vrule width.#2pt
               }
               \hrule height.#2pt
         }
    }
}}

\def\bar{\overline}
\def\real{\mathbb{R}}

\def\lista#1
{{ \itemindent 0.0cm \labelsep .2cm \leftmargin 0.8cm \rightmargin
0.0cm \labelwidth 0.6cm \topsep 0.0mm
\parsep 0.0mm
\itemsep 0.0mm
\begin{list}{}
{ \setlength{\leftmargin}{.8cm} \setlength{\rightmargin}{0.0cm}
\setlength{\parsep}{0.0mm} \setlength{\topsep}{.0mm}
\setlength{\parskip}{.0cm} \setlength{\itemsep}{.0cm} }
{#1}\end{list}} }

\newcounter{theorem}

\begin{document}


\title{\bf Tykhonov Well-posedness of a Heat Transfer Problem with Unilateral Constraints}

\vspace{26mm}
	
\vspace{26mm}
{\author{ {\sc Mircea Sofonea}, Perpignan,\  {\sc Domingo A. Tarzia}, Rosario}

\date {}
\maketitle
\thispagestyle{empty}

\vskip 6mm

 {\small{\it Abstract.} We consider an elliptic boundary value problem with unilateral constraints and subdifferential boundary conditions. The problem describes the  heat transfer in a domain $D\subset\R^d$ and its weak formulation  is in the form of a hemivariational inequality  for the temperature field, denoted by $\cP$. We associate to Problem $\cP$ an optimal control problem, denoted by $\cQ$. Then, using  appropriate Tykhonov triples, governed by a nonlinear operator $G$ and a convex $\wK$, we provide results concerning the well-posedness of problems $\cP$ and $\cQ$.   Our main results are Theorems \ref{t1} and \ref{t2}, together with their corollaries. Their proofs  are based on arguments of compactness, lower semicontinuity and pseudomonotonicity.
Moreover, we consider three relevant perturbations of the heat transfer boundary valued problem which lead to penalty versions of Problem $\cP$, constructed with particular choices of $G$ and $\wK$. We prove that Theorems \ref{t1} and \ref{t2} as well as their corollaries can be applied in the study of these problems, in order to obtain various convergence  results.}

\vskip 6mm


\vskip 6mm


\vskip 14mm

\section{Introduction}\label{s1}
\setcounter{equation}0


\vskip0mm
A large number of mathematical models in Physics, Mechanics and Engineering Sciences are expressed in terms of nonlinear boundary value problems which involve inequalities. Their analysis and optimal control was made the object of many books and papers and, therefore, the literature in the field is extensive.
Here we restrict ourselves to mention the books
\cite{BC, G, HR, HMP, HHNL, KO, Kind-St,NP, Pa,P} and, more recently, \cite{C,MOSBOOK,SMBOOK}.
The results presented in~\cite{BC, C, G, HR, HHNL, KO, Kind-St, Pa} concern
the analysis of various classes of variational inequalities and are based on
arguments of monotonicity and convexity, including properties of the subdifferential
of a convex function.
The results in~\cite{HMP, MOSBOOK,NP,P, SMBOOK} concern
the analysis of hemivariational inequalities and are based on properties
of the subdifferential in the sense of Clarke, defined for locally Lipschitz functions,
which may be nonconvex.  Results in the study of optimal control for
variational and hemivariational inequalities   have been discussed in several works,  including \cite{Ba,BT1,BT2,F,Li,MM1,Mig,MP, NST, ZP,PK,SS,S1, SBT}} and, more recently, in \cite{SBT,ST2}.

The main task in the analysis and control  of inequality problems is to provide their well-posedness. The concepts of well-posedness   vary from problem to problem and  from author to author.  A few examples are the concept of well-posedness in the sense of Hadamard for partial differential equations, the concept of well-posedness in the sense of  Tykhonov for a minimization problem, the concept of well-posedness in the sense of Levitin-Polyak for a constrainted  optimization problem, among others.
The concept of Tykhonov well-posedness (well-posedness, for short)
was introduced in \cite{Ty} for a minimization problem and then
it has been generalized for different optimization problems, see for instance \cite{CKR, DZ, H, HY, L, Z}.  It has been extended in the recent years to various mathematical problems like inequalities, inclusions, fixed point and saddle point problems. Thus, the well-posedness  of variational inequalities
was studied for the first time in \cite{LP1,LP2} and the study of well-posedness  of hemivariational inequalities was initiated in \cite{GM}.  References in the field include \cite{SX8,XHW}, among others. A general framework which unifies the view  on the
well-posedness  for abstract problems in metric spaces  was recently introduced in \cite{SX9}.
Moreover, the  Tykhonov well-posedness of an antiplane shear problem was studied in our recent paper \cite{ST1}.

In this paper we consider a stationary boundary value problem
which describes the  heat transfer in a domain $D\subset\R^d$, $d=1,2,3$. In particular cases, the problem was already considered in \cite{BT1,BT2}. There, the heat flux was assumed to be given on the part $\Gamma_3$ of the boundary of $D$ and, therefore, the weak formulation of the problem was in a form of the variational inequality for the temperature field. In contrast, in this current paper we model the heat transfer on $\Gamma_3$ by using a subdifferential boundary condition, governed  by a locally Lipschitz potential function.
As a consequence, the weak formulation of the   problem leads to an elliptic hemivariational inequality, denoted by $\cP$.  Moreover, we associate to this inequality   an optimal control problem, denoted by $\cQ$.

Our aim in this paper is  twofold. The first one is to study the well-posedness of the hemivariational inequality $\cP$. The second one is to  study the well-posedness of the associated optimal control problem $\cQ$.  To this end, for both problems, we use specific
Tykhonov triples, constructed with approximating sets
govened by a penalty-type operator and a set  of constraints. Proving the well-posedness of the corresponding  hemivariational inequality  and the associated optimal control problem in this functional setting is non standard and represents the  main trait of novelty of our paper. Moreover, the well-posedness of these problems implies general convergence results  which  allows us to deduce the weak-strong dependence of the solution to Problem $\cP$ with respect to the data,  its approch  by a penalty-like method as well as the weak compactness of the set of solution to Problem $\cQ$.  All these results
represent a continuation of our previous results obtained in \cite{BT1,BT2}.

The rest of the manuscript is structured as follows. In Section \ref{s1n} we recall some preliminary material on hemivariational inequalities and well-posedness in the sense of Tykhonov. In Section \ref{s2} we introduce the heat transfer problem, list the assumptions on the data and state its variational formulation,  $\cP$. Then, we  introduce the associate optimal control problem $\cQ$.
In Section \ref{s3}
we prove our first result, Theorem \ref{t1}, which state the well-posedness of Problem $\cP$.
In Section \ref{s4} we  prove our second result, Theorem \ref{t2}, which state
the weakly generalized well-posedness of Problem $\cQ$. Both theorems are completed by several corollaries. Their proofs are based on arguments of pseudomonotonicity, compacteness, and lower semicontinuity.
Finally, in Section \ref{s5} we consider three relevant versions of the heat transfer problem with penalty conditions. We prove that Theorems \ref{t1}, \ref{t2} and their corollaries can be applied in the study of these problems
in order to deduce various convergence results.

\section{Preliminaries}\label{s1n}
\setcounter{equation}0

We start with some notation and preliminaries on hemivariational inequalities and send the reader to \cite{MOSBOOK,NP,SMBOOK} for more details on the material presented below. Let $(X,\|\cdot\|_X)$ be a reflexive Banach space and let $X^*$ denote its dual. We use the symbol $\langle\cdot,\cdot\rangle$ for the duality pairing between $X^*$ and $X$. The limits, upper and lower limits below are considered as $n\to\infty$, even if we do not mention it explicitely. The symbols ``$\rightharpoonup$"  and ``$\to$"
denote the weak and the strong convergence in $X$ or $X^*$.

\medskip
For real valued functions defined on $X$ we recall the following definitions.

\begin{Definition}\label{def}
	{ A function $j \colon X \to \real$ is said to be
		{\rm locally Lipschitz}, \index{function!locally Lipschitz}
		if for every
		$x \in X$, there exists $U_x$ a neighborhood of $x$ and a constant $L_x>0$
		such that\
		$
		|j(y) - j(z)| \le L_x \| y - z \|_X
		$\
		for all $y$, $z \in U_x$. For such functions the
		{\rm generalized (Clarke) directional derivative} of $j$ at the point
		$x \in X$ in the direction $v \in X$ is defined
		by
		\begin{equation*}
		j^{0}(x; v) = \limsup_{y \to x, \ \lambda \downarrow 0}
		\frac{j(y + \lambda v) - j(y)}{\lambda}.
		\end{equation*}
		The {\rm generalized gradient (subdifferential)} of $j$ at $x$
		is a subset of the dual space $X^*$ given by
		\begin{equation*}
		\partial j (x) = \{\, \zeta \in X^* \mid j^{0}(x; v) \ge
		{\langle \zeta, v \rangle} \quad \forall\, v \in X \, \}.
		\end{equation*}
		The function $j$ is said to be {\rm regular} (in the sense
		of Clarke) at the point $x \in X$ if for all $v \in X$ the one-sided directional
		derivative $j' (x; v)$ exists and $j^0(x; v) = j'(x; v)$.}
\end{Definition}

We shall use the following properties of the generalized directional derivative and the generalized gradient.

\begin{Proposition}\label{subdiff}
	Assume that $j \colon X \to \real$ is a locally Lipschitz function. Then the following hold:
	
	\medskip{\rm (a)}
	For every $x \in X$, the function
	$X \ni v \mapsto j^0(x;v) \in \real$ is positively  homogeneous,
	i.e., $j^0(x; \lambda v) = \lambda j^0(x; v)$ for all
	$\lambda \ge 0$, $v\in X$.
	
	\medskip{\rm (b)}
	For every $v \in X$, we have
	$j^0(x; v) = \max \, \{ \, \langle \xi, v \rangle
	\ :\ \xi \in \partial j(x) \, \}$.
\end{Proposition}

We shall use these definitions and properties  both in the case when $X=\R$ and $X=V$ where $V$ is a Sobolev-type space which will defined in Section \ref{s2}. Next, we proceed with the definition of some classes of nonlinear operators.

\begin{Definition}\label{defm}
	An operator $A \colon X \to X^*$ is said to be:
	
	\smallskip
	a)  monotone,
	if for all $u$, $v \in X$, we have $\langle Au - A v, u-v \rangle \ge 0$;
	
	\smallskip
	b) strongly monotone, if there exists $m_A > 0$ such that
	\[\langle Au - Av, u-v \rangle
	\ge
	m_A \| u - v \|_X^{2} \quad \mbox{for all}
	\ u,\, v \in X;\]
	
	\smallskip
	c) bounded, if $A$ maps bounded sets of $X$
	into bounded sets of $X^*$;
	
	\smallskip
	d) pseudomonotone,
	if it is bounded and $u_n \rightharpoonup u$  in $X$ with
	$$\displaystyle \limsup\,\langle A u_n, u_n -u \rangle \le 0$$
	imply\ \ $\displaystyle \liminf\, \langle A u_n, u_n - v \rangle\ge \langle A u, u - v \rangle$\ for all $v \in X$;
	
	\smallskip
	e) demicontinuous,
	if $u_n \to u$ in $X$ implies $A u_n \rightharpoonup Au$ in $X^*$.
	
\end{Definition}

We shall use the following result related to the pseudomonotonicity of operators.

\begin{Proposition}\label{pseudo}
	Assume that the operator $A\colon X \to X^*$
	is bounded, demicontinuous and monotone. Then $A$ is pseudomonotone.
	
\end{Proposition}

We turn now to the study of hemivariational inequalities of the form
\begin{equation}\label{1}
u\in K,\quad\langle A u, v - u \rangle + j^0(u; v - u)
\ge \langle f, v - u \rangle \qquad\forall\,  v \in K.
\end{equation}
and, to this end, we consider the following
hypotheses on the data.
\begin{eqnarray}
&&\label{KK}
K \ \mbox{is nonempty, closed and convex subset of} \ X.
\\[4mm]
&&\label{A}
\left\{\begin{array}{l}A \colon X \to X^* \ \mbox{is  pseudomonotone and}\\
\mbox{strongly monotone with constant}\ m_A>0.
\end{array}\right.
\\[4mm]
&&\label{5jj} \left\{
\begin{array}{l}
j \colon X \to \real \ \mbox{is such that}\\ [2mm]
\ \ {\rm (a)} \
j \ \mbox{is locally Lipschitz.} \\ [2mm]
\ \ {\rm (b)} \
\| \xi \|_{X^*} \le \widetilde{c}_0 + \widetilde{c}_1 \, \| v \|_X \ \mbox{for all}
\ v \in X,\ \xi\in \partial j(v), \ \mbox{with} \ \widetilde{c}_0, \widetilde{c}_1 \ge 0. \\ [2mm]
\ \ {\rm (c)} \
j^0(v_1; v_2 - v_1) + j^0(v_2; v_1 - v_2) \le \alpha_j \, \| v_1 - v_2 \|_X^2 \\ [1mm]
\qquad \ \mbox{for all} \ v_1, v_2 \in X,\ \mbox{with}\ \alpha_j > 0.
\end{array}
\right.
\end{eqnarray}
\vspace{-3mm}
\begin{eqnarray}
&&\label{smal}
\alpha_j < m_A.\\ [2mm]
&&\label{ff}
f \in X^*.
\end{eqnarray}

\medskip


The unique solvability of the variational-hemivariational inequality (\ref{1}) is provided by the following  result proved in \cite{SMBOOK}.

\begin{Theorem}\label{t0}
	Assume  $(\ref{KK})$--$(\ref{ff})$. Then, there exists a unique solution to inequality $(\ref{1})$.

\end{Theorem}

We end this section by recalling  some  preliminaries concerning the concept of well-posedness in the sense of Tykhonov. For more details in the matter we send the  reader to \cite{SX9}.

Assume  that $\cM$ is an abstract mathematic object called generic  ``problem", associated to a metric  space $X$.  Problem $\cM$ could be an equation, a minimization problem, a fixed point problem, an optimal control problem, an inclusion or an inequality problem. Its rigorous statement vary from exemple to example. We  associate to $\cM$ the concept of ``solution" which follows from the context and which
will be clearly defined in each exemple we present below in this paper.
We now introduce the following definitions.

\begin{Definition}\label{def0}
	A Tykhonov triple for Problem $\cM$ is a mathematical object of the form $\cT=(I,\Omega,\cC) $ where  $I$ is a given nonempty set, $\Omega:I\to 2^X-\{\emptyset\}$ and $\cC\subset\cS(I)$, $\cC\ne\emptyset$.	
\end{Definition}

Note that in this definition  and below in this paper $\cS(I)$ represents the set of sequences of elements of $I$ and  $2^X$ denotes the set of parts of the space $X$. A typical element of $I$ will be denoted by $\theta$ and a typical element of $\cS(I)$ will be denoted by $\{\theta_n\}$. We refer to the set $I$ as the set of indices. Moreover, for any $\theta\in I$ we refer to the set $\Omega(\theta)\subset X$ as  an approximating set and $\cC$ will represent the so-called convergence criterion.

\begin{Definition}\label{def1}	
	Given a  Tykhonov triple $\cT=(I,\Omega,\cC)$, 	a sequence $\{u_n\}\subset X$ is called an approximating sequence  if there exists a sequence  $\{\theta_n\}\subset\cC$ such that   	$u_n\in \Omega(\theta_n)$ for each $n\in\mathbb{N}$.
\end{Definition}

Note that approximating sequences always exist since, by assumption, $\cC\ne\emptyset$ and, moreover, for
any sequence $\{\theta_n\}\subset\cC$ and any $n\in\mathbb{N}$, the set $\Omega(\theta_n)$ is not empty. Therefore, Definition \ref{def1} above make sense.

\begin{Definition}\label{def2} 	Given a  Tykhonov triple $\cT=(I,\Omega,\cC)$, Problem $\cM$  is said to be :

\medskip	
a) (strongly) well-posedness  if it has a unique solution  and every approximating sequence  converges  in $X$ to this solution;

\medskip
b) weakly  well-posedness  if it has a unique solution  and every approximating sequence converges weakly in $X$ to this solution;

\medskip
 c) weakly generalized well-posedness  if it has at least one  solution  and every approximating sequence contains a subsequence which converges weakly in $X$ to some point of the solution  set.
\end{Definition}

We remark that  the concept of  approximating sequence above depends on the Tykhonov triple  $\cT$ and, for this reason, we sometimes refer to approximating sequences with respect $\cT$. As a consequence, the concepts of strongly, weakly and weakly generalized well-posedness depend on the Tykhonov triple $\cT$ and, therefore, we  refer to them as strongly, weakly and weakly generalized well-posedness with respect $\cT$, respectively.

\section{The heat transfer problem}\label{s2}
\setcounter{equation}0

In this section we introduce the heat transfer problem, list the assumptions on the data, derive its variational formulation and state our optimal control problem. The problem under consideration is the following.

\medskip\noindent{\bf Problem $\cH$}. {\it
	Find a temperature field $u:D\to\R$ such that}
\begin{eqnarray}
&&\label{d1}u\ge 0,\qquad -\Delta u-f\ge 0,\qquad u(\Delta u+f)=0\qquad{\rm a.e.\ in\ }D,\\ [2mm]
&&\label{d2}u=0\hspace{24mm}{\rm a.e.\ on\ }\Gamma_1,\\ [2mm]
&&\label{d3}u=b\hspace{24mm}{\rm a.e.\ on\ }\Gamma_2,\\ [2mm]
&&\label{d4}-\frac{\partial u}{\partial\nu}\in\partial j_\nu(u)\hspace{8mm}{\rm a.e.\ on\ }\Gamma_3.
\end{eqnarray}

\medskip

In (\ref{d1})--(\ref{d4}) and below $D$ is a bounded domain in $\real^d$ ($d=1,2,3$ in applications) with smooth boundary $\partial D=\Gamma_1\cup\Gamma_2\cup\Gamma_3$ and outer normal unit $\bnu$. We assume that $\Gamma_1$, $\Gamma_2$, $\Gamma_3$ are measurable sets such that $\Gamma_i\cap\Gamma_j=\emptyset$ for $i\ne j$ and, moreover, $meas\,(\Gamma_1)>0$. In addition, in (\ref{d1})--(\ref{d4}) we do not mention the dependence of the different functions on the spatial variable $\bx\in D\cup\partial D$.
The functions $f$, $b$ and $j_\nu$ are given and will be described below. Here we restrict ourselves to mention that  $f$ represents the internal energy, $b$ is a the prescribed temperature  field on $\Gamma_2$,  $\frac{\partial u}{\partial\nu}$ denotes the normal derivative of $u$ on $\Gamma_3$ and  $\partial j_\nu$ denotes the Clarke subdifferential of the potential function $j_\nu$, assumed to be locally Lipschitz.
Note that Problem $\cH$ represents an extension
of the problem considered in \cite{BT1}, where the boundary condition (\ref{d4}) was of the form
\begin{equation}\label{ee}-\frac{\partial u}{\partial\nu}=q \hspace{7mm}{\rm a.e.\ \ on\ \ }\Gamma_3,
\end{equation}
with a given function $q:\Gamma_3\to\real$. It is obvious to see that condition (\ref{ee}) can be recovered by condition (\ref{d4}) with an appropriate choice of $j_\nu$.

\medskip

For the variational analysis of Problem $\cH$ we use standard notation for $C^1$, Lebesque and Sobolev spaces. We use the symbols $``\to"$ and $"\rightharpoonup"$ to indicate the strong and weak convergence in various spaces which will be indicated below. We also use $``\to"$  for the convergence in $\real$. Moreover,  we shall use the space
\[V=\ \{\,v\in H^1(D)\ :\ v=0\quad{\rm on}\quad\Gamma_1\,\},\]
endowed with the inner product  of the space $H^1(D)$, denoted by $(\cdot,\cdot)_V$, and the
associated norm, denoted by $\|\cdot\|_V$. Since $meas (\Gamma_1)>0$, it is well known that $(V,(\cdot,\cdot)_V)$ is a real Hilbert space.
In addition, by the
Friedrichs-Poincar\'e inequality and the Sobolev trace theorem we have
\begin{eqnarray}
&&\label{FP}
\|v\|_{V}\leq c_0\,\|\nabla
v\|_{L^2(D)^d}, \\ [2mm]
&&\label{trace}
\|v\|_{L^2(\Gamma_3)}\le c_3\,\|v\|_V
\end{eqnarray}
for all $v\in V$, respectively. Here and below in this paper
$c_0$ and $c_3$ are positive constants which depend on $D$, $\Gamma_1$ and $\Gamma_3$.
We denote in what follows by $V^*$ the dual of $V$, by $\langle\cdot,\cdot\rangle$ the duality paring between $V^*$ and $V$  and by $0_V$ the zero element of $V$.

\medskip

We now list the assumptions on the data of Problem ${\cal H}$. First,  for the functions $f$ and $b$
we assume that
\begin{eqnarray}
&&\label{f}
f\in L^2(D),\\ [2mm]
&&\label{b}\left\{\begin{array}{ll} b\in L^2(\Gamma_2)\ \mbox{and there exists $u_b\in V$ such that}\\
u_b\ge0  \ \ {\rm in}\ \ D\ {\rm and}\  u_b=b\ \ {\rm on}\ \ \Gamma_2.
\end{array}\right.
\end{eqnarray}
Moreover, for the potential function $j_\nu$ we  assume the following.

\begin{eqnarray}
&&\left\{\begin{array}{ll} j_\nu \colon \Gamma_3 \times \real \to
\real \ \mbox{is such that} \\ [2mm]
\ \ {\rm (a)} \
j_\nu (\cdot, r) \ \mbox{is measurable on} \ \Gamma_3 \ \mbox{for all} \ r \in \real
\ \mbox{and there}\qquad \\
\qquad  \mbox{exists} \ {\bar{e}} \in L^2(\Gamma_3)
\ \mbox{such that} \ j_\nu(\cdot, {\bar{e}} (\cdot)) \in L^1(\Gamma_3). \\ [2mm]
\ \ {\rm (b)} \
j_\nu(\bx, \cdot) \ \mbox{is locally Lipschitz on} \ \real \ \mbox{for a.e.}
\ \bx \in \Gamma_3. \\ [2mm]
\ \ {\rm (c)} \
| \xi | \le {\bar{c}}_0 + {\bar{c}}_1 \, | r |
\ \mbox{for a.e.} \ \bx \in \Gamma_3, \\
\qquad \mbox{ for all} \ r \in \real,\ \xi\in\partial j_\nu(\bx,r), \ \mbox{with}
\ {\bar{c}}_0, \, {\bar{c}}_1 \ge 0. \\ [2mm]
\ \ {\rm (d)} \
j_\nu^0(\bx, r_1; r_2 - r_1) + j_\nu^0(\bx, r_2; r_1 - r_2) \le
\alpha_{j_\nu} \, | r_1 - r_2 |^2 \\
\qquad \mbox{ for a.e.} \ \bx \in \Gamma_3, \ \mbox{all} \ r_1,\,
r_2 \in \real, \ \mbox{with} \ \alpha_{j_\nu} \ge 0,\\[2mm]
\ \ {\rm (e)\  either}\ j_\nu(\bx,\cdot)\ {\rm or}\  -j_\nu(\bx,\cdot)
\ {\rm is \ regular \ on } \ \mathbb{R}\ \ {\rm for\ a.e.}\  \bx\in\Gamma_3.
\end{array}
\right. \label{j}
\end{eqnarray}
Finally, we assume that the following smallness condition holds:
\begin{equation}\label{sm}
\alpha_{j_\nu}c_0^2c_3^2<1.
\end{equation}

Examples of sets $D$, $\Gamma_1$, $\Gamma_2$ and functions $b$ for which such assumption is satisfied can be easily constructed. Below, we restrict ourselves to the following  ones.
\begin{Example} Let $\alpha>0$, $\beta>0$,  and let
	\begin{eqnarray*}
		&&D=\{\,\bx=(x_1,x_2)\in\mathbb{R}^2 \ :\ 0< x_1<\alpha,\ \ 0< x_2<\beta\,\},\\
		&&\Gamma_1=\{\,\bx=(x_1,x_2)\in\mathbb{R}^2 \ : \
		x_1=0,\ \ 0\le x_2\le \beta\,\},\\
		&&\Gamma_2=\{\,\bx=(x_1,x_2)\in\mathbb{R}^2 \ :\ x_1=\alpha,\ \
		0\le x_2\le \beta\,\}.
	\end{eqnarray*}
Assume that $\varphi$ is a positive function with regularity $\varphi\in C^1([0,\beta])$
and consider the functions $b:\Gamma_2\to\mathbb{R}$, $u_b: D\to\mathbb{R}$
defined by $b(x_1,x_2)=\varphi(x_2)$, $u_b(x_1,x_2)=\frac{x_1}{\alpha}\,\varphi(x_2)$.
Then, it is easy to see that condition $(\ref{b})$ is satisfied.
\end{Example}

\begin{Example}
	Let $\alpha>0$, $\beta>0$, $\gamma>0$ and let
	\begin{eqnarray*}
		&&D=\{\,\bx=(x_1,x_2,x_3)\in\mathbb{R}^3 \ :\ 0< x_1<\alpha,\ \ 0< x_2<\beta,\ \ 0< x_3<\gamma\},\\
		&&\Gamma_1=\{\,\bx=(x_1,x_2,x_3)\in\mathbb{R}^3 \ : \
		x_1=0,\ \ 0\le x_2\le \beta,\ \ 0\le x_3\le\gamma\,\},\\
		&&\Gamma_2=\{\,\bx=(x_1,x_2,x_3)\in\mathbb{R}^3 \ :\ x_1=\alpha,\ \
		0\le x_2\le \beta, \ \ 0\le x_3\le\gamma\,\}.
	\end{eqnarray*}
	Assume that $\varphi$ is a positive function with regularity $\varphi\in C^1([0,\beta]\times [0,\gamma])$
	and consider the functions $b: \Gamma_2\to\mathbb{R}$, $u_b:D\to\mathbb{R}$
	defined by $b(x_1,x_2,x_3)=\varphi(x_2,x_3)$, $u(x_1,x_2,x_3)=\frac{x_1}{\alpha}\,\varphi(x_2,x_3)$.
	Then, it is easy to see that condition $(\ref{b})$ is satisfied.
\end{Example}

Note also that various examples of functions $j_\nu$ which satisfy condition (\ref{j}) can be found in \cite{ACTA,MOSBOOK}.

\medskip
Next, we define the set $K\subset V$, the bilinear form  $a:V\times V\to \R$ and the function $j:V\to\R$ by equalities
\begin{eqnarray}
&&\label{K} K=\ \{\,v\in V\ :\ v\ge 0\ \ {\rm in}\ \ D,\quad v=b\ \ {\rm on}\ \ \Gamma_2\,\},
\\ [3mm]
&&\label{a} a(u,v)=\int_{D}\,\nabla u\cdot\nabla
v\,dx\qquad\forall\,u,\,v\in V,\\[2mm]
&&\label{jj} j(v)=\int_{\Gamma_3}j_\nu(v)\,da\qquad\forall\, v\in V.
\end{eqnarray}

It follows from assumption (\ref{b}) that the set $K$ is not empty.  Moreover, Lemma 8 in \cite{SMBOOK} provides the following result.

\begin{Proposition}\label{pr1}
Under assumption $(\ref{j})$ the
function $j:V\to\R$ defined by {\rm (\ref{jj})} satisfies conditions $(\ref{5jj})$ on the space $X=V$ with $\widetilde{c}_0 = \sqrt{2\, {\it meas} (\Gamma_3)} \, {\bar{c}}_0c_3$, $\widetilde{c}_1 = \sqrt{2} \, {\bar{c}}_1  c_3^2$ and $\alpha_j=\alpha_{j_\nu}c_3^2$. In addition,
\begin{equation}
\label{5j}
j^0(u;v)=\displaystyle\int_{\Gamma_3}j_\nu^0(u,v)\,da\qquad\forall\, u,\ v\in V.	
\end{equation}

\end{Proposition}

Next, we turn to the variational formulation of Problem $\cH$. To this end we assume that $u$ is a regular function which satisfies (\ref{d1})--(\ref{d4}) and we consider an arbitrary element  $v\in K$.  Then, using an integration by parts if follows that
\[\int_D\nabla u\cdot\nabla(v-u)\,dx+\int_D\Delta u\cdot(v-u)\,dx=\int_{\partial D}\frac{\partial u}{\partial\nu}\,(v-u)\,da.\]
Therefore, by (\ref{d1})--(\ref{d3}) we find that
\begin{equation}
\label{z2}\int_D\nabla u\cdot\nabla(v-u)\,dx\ge\int_D f(v-u)\,dx+\int_{\Gamma_3}\frac{\partial u}{\partial\nu}\,(v-u)\,da
\end{equation}
and, using the boundary condition (\ref{d4}) combined with definition  of the Clarke subdifferential, we obtain that
\begin{equation}\label{z3}
\int_{\Gamma_3}\frac{\partial u}{\partial\nu}\,(v-u)\,da\ge-\int_{\Gamma_3}j^0_\nu(u;v-u)\,da.
\end{equation}
We now combine inequalities (\ref{z2}) and (\ref{z3}), then use notation (\ref{a}) and equality (\ref{5j}) to see that
\begin{equation*}
a(u,v-u)+j^0(u;v-u)\geq (f,v-u)_{L^2(D)}.
\end{equation*}
Finally, using this inequality and the regularity $u\in K$, guaranteed by (\ref{d1})--(\ref{d3}), we deduce  the following variational formulation of Problem $\cH$.

\medskip\noindent{\bf Problem} $\cP$. {\it
	Find  $u\in V$ such that}
\begin{equation}\label{hv}
u\in K,\quad a(u,v-u)+j^0(u;v-u)\geq (f,v-u)_{L^2(\Omega)}\quad\forall\, v\in K.
\end{equation}

\medskip	
We refer to  inequality (\ref{hv})  as a hemivariational inequality with unilateral constraint.
A function $u\in V$ which satisfies (\ref{hv})
and is called a {\it weak solution} to the heat transfer problem ${\cal H}$.

\medskip
We now introduce the set of admissible pairs for inequality Problem $\cP$ defined by
\begin{equation}\label{ad}
{\cal V}_{ad} = \{\,(u, f)\in K\times L^2(D) \ \mbox{such that}\   (\ref{hv})\  \mbox{holds}\,\}.
\end{equation}
We also assume that $\cL:V\times L^2(D)\to\R$ is a given cost functional which  satisfies the following conditions.
\medskip
\begin{eqnarray}
&&\label{L1}\left\{\begin{array}{l}
\mbox{For all sequences }\{u_n\}\subset V\mbox{ and
}\{f_n\}\subset
L^2(D)\mbox{ such that}\\[2mm]
u_n\rightarrow u\mbox{\ \ in\ \ }V,\ \ f_{n}\rightharpoonup
f \mbox{\ \ in\ \ }
L^2(D),\ \mbox{we have} \\[3mm]
 \displaystyle\liminf_{n\to
	\infty}\,{\cal L}(u_n,f_n)\ge {\cal L}(u,f).

\end{array}\right.
\\[4mm]
&&\label{L2}\left\{ \begin{array}{l}
\mbox{ There exists}\ h: L^2(D)\to\R\ \mbox{such that}
\\[2mm]
{\rm (a)}\ \ {\cal L}(u,f)\ge h(f)\quad \forall\,u\in V,\  f\in L^2(D).\\ [2mm]
{\rm (b)}\ \ \|f_{n}\|_{L^2(D)}\to+\infty\ \Longrightarrow\ h(f_n)\to \infty.
\end{array}\right.
\end{eqnarray}

\begin{Example}\label{e0}
	A typical example of function ${\cal L}$ which satisfies conditions $(\ref{L1})$ and $(\ref{L2})$ is obtained by taking
	\[
	{\cal L}(u,f)=g(u)+h(f)\qquad\forall\, u\in V,\ f\in L^2(D),
	\]
	where $g:V\to\R$ is a continuous positive function and $h: L^2(D)\to \R$ is a weakly lower semicontinuous coercive function, i.e., it satisfies condition $(\ref{L2}){\rm (b)}$.		
\end{Example}

We also assume that $\cL:V\times L^2(D)\to\R$ is a
 given cost function and we associate to Problem $\cP$ the following optimal control problem.

\medskip\medskip\noindent
{\bf Problem} ${\cal Q}$.  {\it Find $(u^*, f^*)\in {\cal V}_{ad}$ such that}
\begin{equation}\label{3d}
\cL(u^*, f^*)=\min_{(u,f)\in {\cal V}_{ad}}\cL(u,f).
\end{equation}

\medskip

Under the previous assumptions, the unique solvability of Problem $\cP$ as well as the solvability of Problem $\cQ$ can be obtained by using standard arguments, already used in \cite{SMBOOK} and \cite{MM1,MM3}, respectively.
Our aim in what follows is to study the well-posedness of these problems  and to derive some consequences. To this end, as already mentionned in Section \ref{s1n}, for each problem we need to prescribe a
Tykhonov triple, based on three ingredients: a set of indices, a familly of approximating sets and a convergence criterion for the sequences of indices.

\section{A well-posedness result}\label{s3}
\setcounter{equation}0

This section is devoted to the well-posedness of the Problem $\cP$, in the sense precised in Section \ref{s1n}, with $X=V$. In order to construct an appropriate Tykhonov triple for this problem,  we consider a set $\wK$, a penalty operator $G:V\to V^*$ and a penalty parameter $\lambda$ such that the following conditions hold.
\begin{eqnarray}
&&\label{wK}
\left\{\begin{array}{ll}
\mbox{\rm (a)}\quad  \wK \ \mbox{\rm is a  closed convex subset of} \ V.\\ [2mm]
\mbox{\rm (b)}\quad  K\subset\wK.
\end{array}\right.\\[4mm]
&&\label{G}
\left\{\begin{array}{ll}
\mbox{\rm (a)}\quad G:V\to V^* \ \ \mbox{is a
bounded, demicontinuous and  monotone operator}.\\ [2mm]
\mbox{\rm (b)}\quad\langle Gu,v-u\rangle\le 0\qquad\forall\, u\in\wK,\ v\in K.\\ [2mm]
\mbox{\rm (c)}\quad u\in \wK,\quad \langle Gu,v-u\rangle=0\quad\forall\,v\in K\ \ \Longrightarrow\ \  u\in K.
\end{array}\right.
 \end{eqnarray}


\medskip
We start with the following result.

\begin{Proposition}\label{pr}
Assume that $(\ref{b})$--$(\ref{sm})$, $(\ref{wK})$ and $(\ref{G})(a)$  hold. Then for each $\lambda>0$ and $f\in L^2(D)$ there exists a unique element $u=u(\lambda,f)$ such that
\begin{equation}\label{hvp}
u\in \wK,\quad
a(u,v-u)+\frac{1}{\lambda}\,\langle Gu,v-u\rangle+j^0(u,;v-u)\geq ({f},v-u)_{L^2(D)}\quad \forall\, v\in \wK.
\end{equation}
\end{Proposition}

\medskip\noindent{\it Proof.}  Define the operator
$A:V\to V^*$ by equality
\begin{equation*}\label{z10}
\langle Au,v\rangle=a(u,v)+\frac{1}{\lambda}\,\langle Gu,v\rangle\qquad\forall\, u,\ v\in V.
\end{equation*}
Then, using the definition (\ref{a}) of the form $a$, inequality (\ref{FP}) and the properties (\ref{G})(a) of the operator $G$,
it is easy to see that the operator $A$ is bounded, demicontinuous and, moreover,
\begin{equation}\label{z11}
\langle Au-Av,u-v\rangle\ge \frac{1}{c_0^2}\,\|u-v\|_V^2\qquad\forall\,u,\, v\in V.
\end{equation}
Therefore, it follows from Proposition \ref{pseudo}  that $A$ is pseudomonotone. In addition, inequality (\ref{z11}) shows that $A$ is strongly monotone with constant $m_A=\frac{1}{c_0^2}$. On the other hand,
the functional $v\mapsto ({f},v)_{L^2(D)}$ is linear
and continuous on $V$ and, therefore, it defines an element in $V^*$. We note that assumption (\ref{b}) implies that $K\ne\emptyset$, hence inclusion  (\ref{wK})(b) guarantees that $\wK\ne\emptyset$. Moreover,
recall assumption
(\ref{wK})(a), Proposition \ref{pr1} and the smallness assumption (\ref{sm}). All these ingradients allow us to apply Theorem \ref{t0}  in order to deduce the unique solvabilty of inequality (\ref{hvp}), which concludes the proof.
\hfill$\Box$

\medskip
We now  take $\cT=(I,\Omega,\cC)$
where
\begin{eqnarray}
\label{I}
&&\hspace{-13mm}I=\{\,\theta=(\lambda,{f})\ :\ \lambda>0,\ {f}\in L^2(\Omega)\,\},\\ [3mm]
&&\label{O}\hspace{-13mm}\Omega(\theta)=\{\, u\in \wK\ {\rm such\ that}\ (\ref{hvp})\ {\rm holds}\,\} \qquad\forall\,\theta=(\lambda,f)\in I,\\ [3mm]
&&\label{C}\hspace{-13mm}\cC=\{\, \{\theta_n\}\ :\   \theta_n=(\lambda_n,f_n)\in I\ \ \forall\, n\in\mathbb{N},\ \
\lambda_n\to 0,\ \
f_{n}\rightharpoonup f\ \ {\rm in}\ L^2(D)\ {\rm  as}\ n\to\infty\,\}.
\end{eqnarray}
Then, using Proposition \ref{pr} we see that for each $\theta=(\lambda,f)$ the set $\Omega(\theta)$  defined by (\ref{O}) is not empty and, therefore, the triple (\ref{I})--(\ref{C}) is a Tykhonov triple in the sense of Definition \ref{def0}.

\medskip
Our main result in this section is the following.

\begin{Theorem}\label{t1} Assume that $(\ref{f})$--$(\ref{sm})$, $(\ref{wK})$ and $(\ref{G})$  hold. Then Problem $\cP$  is  well-posed
with respect to the Tykhonov triple  $(\ref{I})$--$(\ref{C})$.
\end{Theorem}	

\noindent	
{\it Proof.} Following Definition \ref{def2} a),  the proof is carried out in two main steps.

\medskip\noindent i)  {\it Unique solvability of Problem $\cP$.}
Note that assumptions (\ref{wK}) and (\ref{G})(a) are satisfied if $\wK=K$ and $Gv=0_V$ for all $v\in V$, respectively. Moreover, with this particular choice
inequality (\ref{hvp}) reduces to inequality (\ref{hv}), for any $\lambda>0$. Therefore, the existence of a unique solution $u\in K$ to Problem $\cP$ is a direct consequence of Proposition  \ref{pr}.

\medskip\noindent ii) {\it Convergence of approximating sequences.}
To proceed, we
consider  an approximating sequence for the Problem $\cP$, denoted  by $\{u_n\}$. Then, according to Definition \ref{def1}  it follows that
there exists a sequence $\{\theta_n\}$ of elements of $I$, with  $\theta_n=(\lambda_n,f_n)$, such that $u_n\in\Omega(\theta_n)$
for each $n\in\mathbb{N}$ and, moreover,
\begin{eqnarray}
&&\label{C1}\lambda_n\to 0,\\ [2mm]
&&\label{C2}f_{n}\rightharpoonup f\quad {\rm in}\quad L^2(D).
\end{eqnarray}
Note that the inclusion $u_n\in\Omega(\theta_n)$ combined with definition (\ref{O}) imples that
for each $n\in\mathbb{N}$ the following inequality holds:
\begin{eqnarray}
&&\label{hvpp}
u_n\in \wK,\quad
a(u_n,v-u_n)+\frac{1}{\lambda_n}\,\langle Gu_n,v-u_n\rangle+j^0(u_n,;v-u_n)\\[2mm]
&&\qquad\qquad\geq ({f}_n,v-u_n)_{L^2(D)}\quad \forall\, v\in \wK.\nonumber
\end{eqnarray}
Our aim in what follows is to
prove the convergence
\begin{equation}\label{m}
u_n\to u\quad{\rm in}\quad V,\quad{\rm as}\quad n\to\infty.
\end{equation}
To this end we proceed in three intermediate steps that we present below.

\medskip
\noindent {\it
	 {\rm ii-a)} A first weak convergence result}. We claim that there is an element ${\widetilde{u}} \in \wK$ and
a subsequence  of $\{ u_n \}$,
still denoted by $\{ u_n \}$,
such that $u_n \rightharpoonup {\widetilde{u}}$ in $X$, as $n\to \infty$.

\medskip
To prove the claim, we  establish the boundedness of the sequence $\{ u_n \}$ in $V$. Let $n\in\mathbb{N}$ and let $u$ be the solution of Problem $\cP$. We  use assumption (\ref{wK})(b)
and take $v=u$ in (\ref{hvpp}) to see that
\begin{eqnarray*}
	a(u_n, u_n - u)\le
	\frac{1}{\lambda_n}\, \langle Gu_n, u - u_n\rangle
	+j^0(u_n,u-u_n) +
	(f_n, u_n - u)_{L^2(D)}.
\end{eqnarray*}
Then,  using inequalities (\ref{FP}) and (\ref{G})(b)  we obtain that
\begin{eqnarray}
&&\hspace{-5mm}\label{5}
\frac{1}{c_0^2} \, \| u_n  - u\|_V^2 \le
a(u, u - u_n)  +
j^0(u_n;u-u_n)+
(f_n, u_n - u)_{L^2(D)}.
\end{eqnarray}
On the other hand, by  Proposition~\ref{subdiff}(b) we have
\begin{eqnarray*}
&& j^0 (u_n; u - u_n)
= j^0(u_n; u - u_n) + j^0 (u; u_n - u) - j^0 (u; u_n - u)\nonumber \\ [2mm]
&&
\qquad \le j^0(u_n; u - u_n) + j^0 (u; u_n - u)
+ | j^0 (u; u_n - u)| \nonumber \\ [2mm]
&&\quad \qquad
= j^0(u_n; u - u_n) + j^0 (u; u_n - u)  +
|\max\, \{ \,\langle \xi, u_n - u \rangle\, :\, \xi \in \partial j(u)\, \} | \nonumber
\end{eqnarray*}
and, using Proposition \ref{pr1}, we deduce that
\begin{eqnarray}\label{7}
&& j^0 (u_n; u - u_n)
\le
\alpha_{j_\nu}c_3^2 \| u_n - u \|_V^2 + (\widetilde{c}_0 + \widetilde{c}_1 \| u \|_V) \| u_n - u \|_V.
\end{eqnarray}
Finally, note that
\begin{equation}\label{88}
a(u,u-u_n)+(f_n, u_n - u)_{L^2(D)} \le\big(\| u\| _V+ \|f_n \|_{L^2(D)}\big)\| u_n - u\|_V.
\end{equation}
We now combine inequalities (\ref{5})--(\ref{88}) to see that
\begin{eqnarray}
&&\label{9}
\frac{1}{c_0^2} \, \| u_n  - u \|_V^2 \le\big(\| u\| _V+ \|f_n \|_{L^2(D)}\big)\big)\| u_n - u \|_V \\ [2mm]
&&
\quad\qquad +\alpha_{j_\nu}c_3^2 \| u_n - u \|_V^2 + (\widetilde{c}_0 + \widetilde{c}_1 \| u \|_V) \| u_n - u \|_V.\nonumber
\end{eqnarray}

Note that by (\ref{C2}) we know that the sequence $\{f_n\}$ is bounded in $L^2(D)$. Therefore, using inequality (\ref{9}) and
the smallness assumption (\ref{sm}),
we deduce that  there exists a constant $C > 0$
independent of $n$ such that $\| u_n-u\|_X \le C$. This implies that  the sequence $\{u_n\}$ is bounded  in $V$.
Thus, from the reflexivity of $V$, by passing to
a subsequence, if necessary,  we deduce that
\begin{equation}\label{10}
u_n \rightharpoonup {\widetilde{u}} \  \ {\rm in} \ \ V, \ \ \mbox{as}\ \ n\to \infty,
\end{equation}
with some ${\widetilde{u}} \in V$. Moreover, assumption (\ref{wK})(a) and the convergence (\ref{10}) imply that $\wu\in\wK$, which  completes the proof of the claim.

\medskip\noindent {\rm ii-b)} {\it A property of the weak limit.}
Next, we show that ${\widetilde{u}}$ is a solution to
Problem ${\cal P}$.

\medskip
Let $v$ be a given element in   $\wK$ and let $n\in\mathbb{N}$. We use   (\ref{hvpp}) to obtain that
\begin{eqnarray}
&&\label{11b}
\frac{1}{\lambda_n}\, \langle Gu_n, u_n - v\rangle\le a(u_n, v-u_n)+j^0(u_n,v-u_n)+(f_n, u_n - v)_{L^2(D)}.
\end{eqnarray}
Then, using arguments similar to those used in the proof of (\ref{7}), (\ref{88}) and  the boundedness of the sequence  $ \{u_n\}$, we deduce that each term  in the right hand side of inequality (\ref{11b}) is bounded. This implies that there  exists a constant
$M_0>0$  which does not depend on $n$ such that
\begin{equation*}
\langle Gu_n, u_n - v\rangle\le \lambda_nM_0.
\end{equation*}
We now pass to the upper limit in this inequality and use the convergence (\ref{C1}) to deduce that
\begin{equation}\label{13}
\limsup\,\langle Gu_n, u_n -v\rangle\le 0.
\end{equation}
We now take $v=\wu$ in (\ref{13}) and find that
\begin{equation}\label{13n}
\limsup\,\langle Gu_n, u_n -\wu\rangle \le 0.
\end{equation}
Therefore, using  the pseudomonotonicity of the operator $G$, guaranted by assumption (\ref{G})(a) and Proposition \ref{pseudo},  we obtain that
\begin{equation}\label{13nn}
\liminf\,\langle Gu_n, u_n -v\rangle \ge \langle G\wu,\wu-v\rangle.
\end{equation}
We now combine inequalities (\ref{13nn})  and (\ref{13}) to find that
$\langle G {\widetilde{u}}, {\widetilde{u}} - v\rangle\le 0$.  On the other hand, by (\ref{G})(b) we deduce that $\langle G {\widetilde{u}}, {\widetilde{u}} - v\rangle\ge 0$. We conclude from above that \[\langle G {\widetilde{u}}, {\widetilde{u}} - v\rangle=0\] and  we recall that this equality holds for all $v\in \widetilde{K}$. Then,
using assumptions (\ref{wK})(b), (\ref{G})(c), we find that ${\widetilde{u}} \in K$.

Next, we use (\ref{hvpp}), again,    to obtain
that
\begin{eqnarray*}
	&&a(u_n, u_n-v)\le
	\frac{1}{\lambda_n}\, \langle Gu_n, v - u_n\rangle
	+ j^0(u_n,v-u_n)
	+ (f_n, u_n - v)_{L^2(D)}.
\end{eqnarray*}
Therefore, using
assumption (\ref{G})(b)
we find that
\begin{equation}\label{17}
a(u_n, u_n - v) \le
j^0(u_n,v-u_n) + (f_n,u_n -v)_{L^2(D)}.
\end{equation}
Moreover, from  (\ref{10}), (\ref{5j}), the compactness of the trace and the properties of the integral it follows that
\begin{equation}\label{19m}
\limsup j^0(u_n; v - u_n) \le j^0(\widetilde{u},v-\widetilde{u})
\end{equation}
and, in addition,
\begin{equation}\label{19mn}
(f_n,u_n-v)_{L^2(D)}\to(f,\wu-v)_{L^2(D)} .
\end{equation}
We now gather relations (\ref{17})--(\ref{19mn}) to see that
\begin{equation}
\label{21}
\limsup\,a(u_n, u_n - v) \le j^0({\widetilde{u}},v-{\widetilde{u}})
+ (f,\widetilde{u} - v)_{L^2(D)}.
\end{equation}

Now, taking $v= {\widetilde{u}} \in K$ in (\ref{21})  and using Proposition \ref{subdiff}(a)
we obtain that
\begin{equation}\label{24n}
\limsup\, a(u_n, u_n - {\widetilde{u}}) \le 0.
\end{equation}
On the other hand,  using the properties of the form $a$ we have that
\begin{equation}\label{z3n}
a(u_n,v)\to a(\widetilde{u},v)\quad{\rm as}\quad n\to\infty
\end{equation}
and, since\ $a(u_n-\widetilde{u},u_n-\widetilde{u})\ge 0$, we deduce that \[a(u_n,u_n)\ge a(\widetilde{u},u_n)+a(u_n,\widetilde{u})-a(\widetilde{u},\widetilde{u}),\]
for each $n\in\mathbb{N}$.
Using now (\ref{z3n})  and the convergence (\ref{10}) we find that
\begin{equation}\label{z4}
\liminf_{n\to\infty}\, a(u_n,u_n)\ge a(\widetilde{u},\widetilde{u}).
\end{equation}

\noindent
Therefore, combining (\ref{z3n}), (\ref{z4})  we obtain that
\begin{equation}\label{rr}
a({\widetilde{u}}, {\widetilde{u}} - v)\le\liminf\, a(u_n, u_n -v)
\end{equation}
and, using (\ref{rr})  and (\ref{21}) we have
\begin{equation*}
a({\widetilde{u}}, {\widetilde{u}} - v)\\ [2mm]
\le j^0(\widetilde{u},v-\widetilde{u}) +
(f, \widetilde{u} - v)_{L^2(D)}.
\end{equation*}
It follows from here that ${\widetilde{u}} \in K$ is a solution to Problem~${\cal P}$, as claimed.

\medskip\noindent
{\it
	{\rm ii-c)}  A second weak convergence result}. We now prove the weak convergence of the whole sequence $\{u_n\}$.

\medskip
Since Problem~${\cal P}$ has a unique solution $u \in K$, we deduce from the previous step that ${\widetilde{u}} = u$. Moreover, a careful analysis of the proof in step {\rm ii-b)} reveals that every subsequence of $\{ u_n \}$
which converges weakly in $V$ has the  weak limit $u$. In addition, we recall that the sequence  $\{ u_n \}$ is bounded in $V$. Therefore, using  a standard argument we deduce that  the whole sequence $\{ u_n \}$ converges weakly in $V$ to $u$,  as $n\to \infty$, i.e.,
\begin{equation}\label{9n}
{u}_n\rightharpoonup{u}\quad \mbox{ \rm{in} }\
V\quad \mbox{ \rm{as} }\ n\to\infty.
\end{equation}

\medskip\noindent
{\it {\rm ii-d)}  Strong convergence. } In the final step of the proof, we prove that
$u_n \to u$ in $V$, as $n\to\infty$.

\medskip
We take
$v = {\widetilde{u}} \in K$ in (\ref{rr}) and use (\ref{24n}) to obtain
\[
0 \le \liminf\, a(u_n, u_n - {\widetilde{u}})
\le
\limsup\,a(u_n, u_n - {\widetilde{u}}) \le 0,
\]
which shows that $a(u_n, u_n - {\widetilde{u}}) \to 0$,
as $n\to\infty$.
Therefore, using equality $\widetilde{u}=u$, the coercivity of $a$
and the convergence $u_n \rightharpoonup u$  in $V$,
we have
\[
\frac{1}{c_0^2}\, \| u_n - u \|_V^2 \le
a(u_n - u, u_n - u) =
a(u_n, u_n - u)- a( u, u_n - u) \to 0,
\]
as $n\to \infty$. Hence, it follows that $u_n \to u$ in $V$, which completes the proof.
\hfill$\Box$

\medskip\medskip
We now reformulate the convergence result obtained  above.  To this end, we denote by $u(\lambda_n,f_n)$ and $u(f)$ the solution of the hemivariational inequality (\ref{hvpp}) and (\ref{hv}), respectively. Then, if
(\ref{C1})  and (\ref{C2})  hold, the step ii) of the proof of Theorem \ref{t2} implies that
\begin{equation}\label{z40}
u(\lambda_n,f_n)\to u(f)\quad{\rm in}\quad V.
\end{equation}
We shall use this convergence result in various places below in this paper.

\medskip	
We end this section with the following consequences of Theorem \ref{t1}.

 \begin{Corollary}\label{cor1} Assume that $(\ref{b})$--$(\ref{sm})$, $(\ref{wK})$ and $(\ref{G})$  hold. Then the operator $f\mapsto u$ which associates to every $f\in L^2(D)$  the solution $u=u(f)\in V$ of  inequality $(\ref{hv})$ is weakly-strongly continuous, i.e.,
 \begin{equation}\label{z20}
 f_n\rightharpoonup f\quad{\rm in}\quad L^2(D)\quad\Longrightarrow\quad u(f_n)\to u(f)\quad{\rm in}\quad V.
 \end{equation}	
 \end{Corollary}	

 \noindent {\it Proof.}  Let $\{f_n\}\subset L^2(D)$ be a sequence such that $f_n\rightharpoonup f\quad{\rm in}\quad L^2(D)$.
 We take $\wK=K$, $Gv=0_V$ for all $v\in V$ and let $\{\lambda_n\}\subset\R$ be a sequence such that $\lambda_n>0$ for each $n\in\mathbb{N}$ and $\lambda_n\to 0$. Then, it follows that $u(\lambda_n,f_n)=u(f_n)$ for each $n\in\mathbb{N}$ and,  therefore, the convergence (\ref{z20}) is  a direct consequence of the convergence (\ref{z40}).
 \hfill$\Box$

\begin{Corollary}\label{cor2} Assume that $(\ref{b})$--$(\ref{sm})$, $(\ref{wK})$ and $(\ref{G})$  hold. Then, for any
$f\in L^2(D)$  the solution $u(\lambda,f)\in V$ of  inequality $(\ref{hvp})$ converges to the solution
 $u(f$) of  inequality $(\ref{hv})$ as $\lambda\to 0$, that is
	\begin{equation}\label{z21}
	\lambda_n>0,\quad \lambda_n\to 0\quad\Longrightarrow\quad u(\lambda_n,f)\to u(f)\quad{\rm in}\quad V.
	\end{equation}	
\end{Corollary}	

\noindent {\it Proof.}  The convergence (\ref{z21}) is a direct consequence of the convergence (\ref{z40}).
\hfill$\Box$

\medskip
Note that in the case when $\wK=V$ inequality (\ref{hvp})  represents a hemivariational inequality without constraints. In this particular case Corollary \ref{cor2} can be recovered by using convergence results for penalty variational-hemivariational inequalities obtained in \cite{SMBOOK,SX7}. Nevertheless, we stress that the functional framework used in the above mentionned papers does not allow the use of these results in the case when $\wK\ne K$.

\section{Weakly generalized well-posedness results }\label{s4}
\setcounter{equation}0

This section is devoted to the weakly generalized well-posedness of the Problem $\cQ$
in the sense precised in Section \ref{s1n}, with $X=V\times L^2(D)$. We shall construct a relevant example of Tykhonov triple, which will play a crucial role in the study of Problem $\cQ$. To this end, everywhere in this section we assume that $(\ref{b})$--$(\ref{sm})$, $(\ref{L1})$, $(\ref{L2})$, $(\ref{wK})$ and $(\ref{G})$ hold. Then, for any $\lambda>0$ and $f\in L^2(\Omega)$ we consider the hemivariatinal inequality (\ref{hvp}), which represents a perturbation of the hemivariational inequality (\ref{hv}). Therefore, by analogy with (\ref{ad}), we
introduce the set of perturbed admissible pairs for inequality Problem $\cP$, defined by
\begin{equation}\label{adl}
{\cal V}_{ad}^\lambda = \{\,(u, f)\in \wK\times L^2(D) \ \mbox{such that}\   (\ref{hvp})\  \mbox{holds}\,\}.
\end{equation}
Moreover, for each $\mu\ge0$ let $\cL_\mu:V\times L^2(D)\to\R$ be a perturbation of the cost functional, assumed to satisfy conditions (\ref{L1}) and (\ref{L2}), i.e.,
\begin{eqnarray}
&&\label{L1m}\left\{\begin{array}{l}
\mbox{For all sequences }\{u_n\}\subset V\mbox{ and
}\{f_n\}\subset
L^2(D)\mbox{ such that}\hskip24mm\\[2mm]
u_n\rightarrow u\mbox{\ \ in\ \ }V,\ \ f_{n}\rightharpoonup
f \mbox{\ \ in\ \ }
L^2(D),\ \mbox{we have} \\[3mm]
\displaystyle\liminf_{n\to
	\infty}\,{\cal L}_\mu(u_n,f_n)\ge {\cal L}_\mu(u,f).
\end{array}\right.
\\[4mm]
&&\label{L2m}\left\{ \begin{array}{l}
{\rm (a)}\ \ {\cal L}_\mu(u,f)\ge h(f)\quad \forall\,u\in V,\  f\in L^2(D).\\ [2mm]
{\rm (b)}\ \ \|f_{n}\|_{L^2(D)}\to+\infty\ \Longrightarrow\ h(f_n)\to \infty.
\end{array}\right.
\end{eqnarray}

In addition, we assume that the following properties hold.
\begin{eqnarray}
&&\label{L3}\left\{\begin{array}{l}
\mbox{For all sequences }\{u_n\}\subset V,\ \
\{f_n\}\subset
L^2(D) \mbox{ and
}\{\mu_n\}\subset
\R_+\mbox{ such that}\\[2mm]
u_n\rightarrow u\mbox{\ \ in\ \ }V,\ \ f_{n}\rightharpoonup
f \mbox{\ \ in\ \ }
L^2(D),\ \mu_n\to0\ \mbox{we have} \\[3mm]
\displaystyle\liminf_{n\to
	\infty}\,{\cal L}_{\mu_n}(u_n,f_n)\ge {\cal L}(u,f).
\end{array}\right.\\[4mm]
&&\label{L4}\left\{\begin{array}{l}
\mbox{For all sequences }\{u_n\}\subset V\mbox{ and
}\{\mu_n\}\subset
\R_+\mbox{ such that}\\[2mm]
u_n\rightarrow u\mbox{\ \ in\ \ }V,\ \mu_n\to0\ \ \mbox{\rm and for any}\ f\in L^2(D)\  \mbox{we have}\\ [2mm]  \displaystyle\lim_{n\to
	\infty}\, {\cal L}_{\mu_n}(u_n,f)={\cal L}(u,f).
\end{array}\right.
\end{eqnarray}

An example of functionals ${\cal L}_\mu$ and ${\cal L}$
which satisfy conditions (\ref{L3}) and (\ref{L4}) will be provided in the next section. Here we restrict ourselves to remark that
condition (\ref{L4}) implies that for $\mu=0$ the functionals ${\cal L}_\mu$ and
${\cal L}$
are the same, i.e., \begin{equation}\label{id}
{\cal L}_0(u,f)={\cal L}(u,f)\qquad\forall\, u\in V,\ f\in L^2(D).
\end{equation}
It follows from here that if assuming conditions  (\ref{L4}), (\ref{L1m}) and  (\ref{L2m}) for each $\mu\ge 0$ implies that conditions (\ref{L1}) and  (\ref{L2})
hold, too.

\medskip
We start with the following result.

\begin{Proposition}\label{pro}
	Assume  $(\ref{b})$--$(\ref{sm})$, $(\ref{wK})$ and $(\ref{G})$. Then, for each $\lambda>0$ and each $\mu\ge 0$ such that  $(\ref{L1m})$ and  $(\ref{L2m})$ hold,  there exists   a pair $(u^*,f^*)$ such that
	\begin{equation}\label{3dn}
	(u^*,f^*)\in{\cal V}^\lambda_{ad},\qquad \cL_\mu(u^*, f^*)=\min_{(u,f)\in {\cal V}^\lambda_{ad}}\cL_\mu(u,f).
	\end{equation}

\end{Proposition}

\medskip\noindent{\it Proof.} We use standard arguments that we resume in the following. Let  $\lambda>0$, $\mu\ge 0$, let
\begin{equation}\label{omega}
\omega=\inf_{(u,f)\in {\cal V}^\lambda_{ad}} {\cal L}_\mu(u, f)\in\mathbb{R}
\end{equation}
and let $\{(u_n,f_{n})\}\subset {\cal V}^\lambda_{ad}$ be a minimizing sequence for the functional  ${\cal L}_\mu$, i.e.,
\begin{equation}\label{n}
\lim_{n\to\infty}\,{\cal L}_\mu(u_n, f_{n})=\omega.
\end{equation}

We claim that the sequence $\{f_{n}\}$ is bounded in $L^2(D)$.
Arguing by contradiction, assume  that $\{f_{n}\}$ is not bounded in $L^2(D)$. Then, passing to a subsequence still denoted by $\{f_{n}\}$, we have
\begin{equation}\label{8}
\|f_{n}\|_{L^2(D)}\to +\infty \quad\text{as}\quad n\to +\infty
\end{equation}and, using  (\ref{L2m}) it turns that
\begin{equation}\label{mm}
\lim_{n\to +\infty}{\cal L}_\mu(u_n, f_n)= +\infty.
\end{equation}
Equalities (\ref{n}) and (\ref{mm}) imply that $\omega=+\infty$ which is in contradiction with (\ref{omega}).

We conclude from above that the sequence $\{f_{n}\}$ is bounded in $L^2(D)$, as clamed. Therefore, we deduce that there exists $f^*\in L^2(D)$ such that, passing to a subsequence still denoted by $\{f_{n}\}$, we have
\begin{equation}\label{ww}
f_{n}\rightharpoonup f^*\quad\text{in}\quad L^2(D)\quad\text{as}\quad n\to +\infty.
\end{equation}

Let $u^*$ be the solution of the hemivariational inequality (\ref{hvp}) for $f=f^*$.
Then,  by the definition (\ref{adl}) of the set ${\cal V}^\lambda_{ad}$ we have
\begin{equation}\label{sol-opt}
(u^*, f^*)\in {\cal V}^\lambda_{ad}.
\end{equation}
Moreover, using  the convergence (\ref{ww}) and Corollary \ref{cor1}  it follows that
\begin{equation}\label{www}
u_n \to u^*\quad\text{in}\quad V\quad \text{as}\quad n\to +\infty.
\end{equation}
We now use the convergences (\ref{ww}), (\ref{www}) and assumption (\ref{L1m}), to deduce that
\begin{equation}\label{xxx}
\liminf\limits_{n\to +\infty}{\cal L}_\mu(u_n,f_{n})\geq{\cal L}_\mu(u^*,f^*).
\end{equation}
It follows from (\ref{n}) and (\ref{xxx}) that
\begin{equation}\label{xy}
\omega\geq{\cal L}_\mu(u^*,f^*).
\end{equation}
In addition, (\ref{sol-opt})  and (\ref{omega}) yield
\begin{equation}\label{xz}
\omega\le{\cal L}_\mu(u^*,f^*).
\end{equation}
We combine inequalities  (\ref{xy}), (\ref{xz}) then we use (\ref{sol-opt}) see that (\ref{3dn}) holds, which
concludes the proof.
\hfill$\Box$

\medskip
We now take $\cT=(I,\Omega,\cC)$
where
\begin{eqnarray}
\label{Ic}
&&\hspace{-13mm}I=\{\,\theta=(\lambda,\mu),\ :\ \lambda>0,\ \mu\ge 0\,\},\\ [3mm]
&&\label{Oc}\hspace{-13mm}\Omega(\theta)=\{\, (u^*,f^*)\in V\times L^2(D)\ {\rm such\ that}\ (\ref{3d})\ {\rm holds}\,\} \qquad\forall\,\theta=(\lambda,\mu)\in I,\\ [3mm]
&&\label{Cc}\hspace{-13mm}\cC=\{\, \{\theta_n\}\ :\   \theta_n=(\lambda_n,\mu_n)\in I\ \ \forall\, n\in\mathbb{N},\ \
\lambda_n\to 0,\ \
\mu_n\to 0\ \ \ {\rm  as}\ n\to\infty\,\}.
\end{eqnarray}
Then, using Proposition \ref{pro} we see that for each $\theta=(\lambda,\mu)$ the set $\Omega(\theta)$  defined defined by (\ref{3d}) is not empty and, therefore, the triple (\ref{Ic})--(\ref{Cc}) is a Tykhonov triple in the sense of Definition \ref{def0}.

\medskip
Our main result in this section is the following.

\begin{Theorem}\label{t2} Assume
	$(\ref{b})$--$(\ref{sm})$, $(\ref{wK})$, $(\ref{G})$, $(\ref{L3})$, $(\ref{L4})$ and, for all $\mu\ge 0$, assume that $(\ref{L1m})$, $(\ref{L2m})$ hold. Then Problem $\cQ$  is weakly  generalized  well-posed with respect to the Tykhonov triple  $(\ref{Ic})$--$(\ref{Cc})$.
\end{Theorem}

\noindent	
{\it Proof.}  Following Definition \ref{def2} the proof is carried out in two main steps.

\medskip\noindent i) {\it Solvability of Problem $\cQ$.} Recall that for $\wK=K$ and $Gv=0_V$, for all $v\in V$
inequality (\ref{hvp}) reduces to inequality (\ref{hv}). Moreover, (\ref{id}) shows that for $\mu=0$ we have ${\cal L}_\mu={\cal L}$. Therefore, the existence of at least one solution  to Problem $\cQ$ is a direct consequence of Proposition  \ref{pro}.

\medskip\noindent ii) {\it Convergence of approximating sequences.}
To proceed, we
consider  an approximating sequence for the Problem $\cQ$, denoted  by $\{(u_n^*, f_{n}^*)\}$. Then, according to Definition \ref{def1}  it follows that
there exists a sequence $\{\theta_n\}$ of elements of $I$, with $\theta_n=(\lambda_n,\mu_n)$,  such that $(u_n,f_{n})\in\Omega(\theta_n)$
for each $n\in\mathbb{N}$, (\ref{C1}) holds and, moreover,
\begin{equation}\label{C3}
\mu_n\to 0.
\end{equation}
In order to simplify the notation, for each $n\in\mathbb{N}$ we write ${\cal V}^n_{ad}$ and
${\cal L}_n$ instead of  ${\cal V}^{\lambda_n}_{ad}$ and ${\cal L}_{\mu_n}$, respectively.  Then, exploiting the definition (\ref{Oc}) we deduce that  $(u^*_n,f_{n}^*)\in {\cal V}_{ad}^n$ and
\begin{eqnarray}
&&\label{zc}
{\cal L}_n(u^*_n,f^*_{n})\le {\cal L}_n(u_n,f_{n})
\end{eqnarray}
for each couple of functions $(u_n,f_{n})\in  {\cal V}_{ad}^n$, i.e., for each couple of functions $(u_n,f_{n})\in V\times L^2(D)$ which satisfies inequality (\ref{hvp}) in which $f$ is replaced by $f_n$, for each $n\in\mathbb{N}$.

We shall prove that there exists a subsequence of the sequence $\{(u_{n}^{*}, f_{n}^{*})\}$, again denoted by $\{(u_{n}^{*}, f_{n}^{*})\}$, and an element   $(u^*,f^*)\in V\times L^2(D)$ such that
\begin{eqnarray}
&&\label{se1}
f_{n}^* \rightharpoonup f^*\quad\mbox{\rm in}\quad L^2(D)\quad\mbox{\rm as}\quad n\to \infty, \\[2mm]
&&\label{se2}
u_n^*\rightarrow u^*\quad \mbox{\rm in}\quad V\quad\mbox{\rm as}\quad n\to \infty, \\[2mm]
&&\label{se3}
(u^*,f^*)\quad {\rm is\ a\ solution\ of\ Problem}\ {\cal Q}.
\end{eqnarray}
To this end we proceed in four intermediate steps that we present below.

\medskip\noindent {\rm ii-a)} {\it A boundedness result.} We claim that the sequence  $\{f_{n}^*\}$ is bounded in $L^2(D)$.
Arguing by contradiction, assume  that $\{f_n^*\}$ is not bounded in $L^2(D)$. Then, passing to a subsequence still denoted by $\{f_n^*\}$, we have
\begin{equation}\label{8n}
\|f_{n}^*\|_{L^2(D)}\to +\infty\quad\text{as}\quad n\to +\infty.
\end{equation}
We  use assumption (\ref{L2})$_\mu$ of the cost function to deduce that
\begin{equation}\label{20}
\lim_{n\to \infty}{\cal L}_n(u_n^*, f_{n}^*)= +\infty.
\end{equation}
Next, let $f\in L^2(D)$ be given and let ${\wu}_n$ be the solution of the hemivariational inequality
\begin{equation}
\widetilde{u}_n\in \wK\ :\   \label{hvpn}
a(\widetilde{u}_n,v-\widetilde{u}_n)+\frac{1}{\lambda_n}\,\langle G\widetilde{u}_n,v-\widetilde{u}_n\rangle+j^0(\widetilde{u}_n,;v-\widetilde{u}_n)\geq (f,v-\widetilde{u}_n)_V\quad \forall\, v\in \wK,
\end{equation}
i.e., $\widetilde{u}_n=u(\lambda_n,f)$. Note that Proposition \ref{pr} guarantees that this solution exists and is unique, for each $n\in\mathbb{N}$. Moreover,  using (\ref{C1}) and (\ref{z40}) it follows that
\begin{equation}
\widetilde{u}_n\rightarrow u\quad \mbox{\rm in}\quad V\quad\mbox{\rm as}\quad n\to \infty
\end{equation}
where $u=u(f)$
and, by (\ref{C3}) and assumption  (\ref{L4}), we deduce that
\begin{equation}\label{e5}
{\cal L}_n(\widetilde{u}_n, f)\to{\cal L}(u, f) \quad\mbox{\rm as}\quad n\to \infty
\end{equation}

On the other hand (\ref{hvpn})  implies that the pair $(\widetilde{u}_n,f)$ satisfies inequality  (\ref{hvp}), i.e., $(\widetilde{u}_n,f)\in {\cal V}^n_{ad}$. Therefore, (\ref{zc}) implies that
\begin{eqnarray}
&&\label{e6}
{\cal L}_n(u^*_n,f^*_{n})\le {\cal L}_n(\widetilde{u}_n,f)\qquad\forall\, n\in\mathbb{N}.
\end{eqnarray}
We now pass to the limit in (\ref{e6}) as $n\to\infty$ and use  (\ref{20}) and
(\ref{e5}) to obtain a contradiction, which proves the claim.

\medskip\noindent {\rm ii-b)} {\it Two convergence results.} We conclude from step ii-a) that
the sequence $\{f_{n}^*\}$ is bounded in $L^2(D)$ and, therefore,  we can find	a subsequence again denoted by $\{f_{n}^*\}$ and an element $f^*\in L^2(D)$ such that (\ref{se1}) holds.
Denote by $u^*$ the solution of inequality (\ref{hv}) for $f=f^*$ and note that definition (\ref{ad}) implies that
\begin{equation}\label{se6}
(u^*, f^*)\in {\cal V}_{ad}.
\end{equation}
Moreover, since $u_n^*=u(\lambda_n,f_n)$, $u^*=u(f)$, assumption (\ref{C1})  and convergences (\ref{se1}), (\ref{z20}) imply that (\ref{se2}) holds, too.

\medskip\noindent {\rm ii-c)} {\it Optimality of the limit.}
We now prove that $(u^*,f^*)$ is a solution to the optimal control problem ${\cal Q}$.
To this end we use the convergences (\ref{se1}), (\ref{se2}), (\ref{C3}) and assumption (\ref{L3}) to see that
\begin{equation}\label{se8}
{\cal L}(u^*,f^*)\leq\liminf_{n\rightarrow\infty}{\cal L}_n(u_n^*,f_{n}^*).
\end{equation}
Next, we fix a  solution $({u}^*_0,{f}^*_{0})$ of Problem ${\cal Q}$ and,
in addition, for each $n\in\mathbb{N}$  we denote by ${u}_n^0$ the unique element of $\wK$ which satisfies the inequality (\ref{hvp}) with $\lambda=\lambda_n$ and $f=f_0^*$, i.e., $u_n^*=u(\lambda_n,f_0^*)$.
Therefore, $({u}_n^0,{f}^*_{0})\in {\cal V}_{ad}^n$ and, using the optimality of the pair
$(u_n^*,f_{n}^*)$, (\ref{zc}), we find that
\begin{equation*}
{\cal L}_n(u_n^*,f_{n}^*)\leq{\cal L}_n({u}_n^0,{f}^*_{0})\qquad\forall\, n\in\mathbb{N}.
\end{equation*}
We pass to the upper limit in this inequality to see that
\begin{equation}\label{se9}
\limsup_{n\rightarrow\infty}{\cal L}_n(u_n^*, f_{n}^*)\leq \limsup_{n\rightarrow\infty}{\cal L}_n({u}_n^0,{f}^*_{0}).
\end{equation}

Next, since $\lambda_n\to 0$ and $u_0^*=u(f_0^*)$, it follows from (\ref{z40}) that
\begin{equation*}
{u}_n^0 \to {u}^*_0\quad\text{in}\quad V\quad\text{as}\quad n\to\infty.
\end{equation*}
Using now this convergence and assumption (\ref{L4}) yields
\begin{equation}\label{se10}
\lim_{n\rightarrow\infty}{\cal L}_n({u}_n^0,{f}^*_{0})={\cal L}({u}^*_0,f^*_{0}).
\end{equation}
We now use  (\ref{se8})--(\ref{se10}) to see that
\begin{equation}\label{se10n}
{\cal L}(u^*,f^*)\leq {\cal L}({u}^*_0,{f}^*_{0}).
\end{equation}

On the other hand, since $({u}^*_0,{f}^*_{0})$ is a solution of Problem ${\cal Q}$,  we have
\begin{equation}\label{3p}
{\cal L}({u}^*_0,{f}^*_{0})=\min_{(u,f)\in {\cal V}_{ad}} {\cal L}(u,f)
\end{equation} and, therefore,
inclusion (\ref{se6})
implies that
\begin{equation}\label{se11}
{\cal L}({u}^*_0,{f}^*_{0})\le {\cal L}(u^*,f^*).
\end{equation}
We now combine the inequalities  (\ref{se10n}) and (\ref{se11}) to see that
\begin{equation}\label{se16}
{\cal L}(u^*, f^*)={\cal L}({u}^*_0,{f}^*_{0}).
\end{equation}
Next, we use relations (\ref{se6}), (\ref{se16}) and (\ref{3p}) to see that (\ref{se3}) holds.

\medskip\noindent {\rm ii-d)} {\it End of proof.} We remark that the convergences (\ref{se1})  and (\ref{se2})  imply the weak convergence (in the product Hilbert space $V\times L^2(D)$) of the sequence $\{(u_n^*,f_{n}^*)\}$ to the element $(u^*,f^*)$. Theorem \ref{t2} is now a direct consequence of Definition \ref{def2} c). \hfill$\Box$

\medskip
We turn now to some direct consequence of Theorem \ref{t2}.

\begin{Corollary}	\label{cor3}
Assume that $(\ref{b})$--$(\ref{sm})$,
$(\ref{L1})$ and $(\ref{L2})$  hold.  Then, the set of solutions of Problem $\cQ$  is weakly  sequentially compact.		
\end{Corollary}

\noindent	
{\it Proof.} Assume that $\{(u_n^*, f_{n}^*)\}$ is a sequence of solutions to Problem $\cQ$.  Let $\{\lambda_n\}\subset\R$ be such that $\lambda_n>0$ for each $n\in\mathbb{N}$, $\lambda_n\to 0$, and let $\mu_n=0$ for each $n\in\mathbb{N}$.
Also, consider  the particular case when $\wK=K$, $Gv=0_V$ for any $v\in V$, and note that in this case we have
${\cal V}_{ad}^n={\cal V}_{ad}$,
 ${\cal L}_n={\cal L}$,  for each $n\in\mathbb{N}$. It follows from above that
$\{(u_n^*, f_{n}^*)\}$ is an approximating sequence for Problem $\cQ$. Therefore, the step (ii) in the proof of Theorem \ref{t2} implies that there exists a subsequence of the sequence $\{(u_{n}^{*}, f_{n}^{*})\}$, again denoted by $\{(u_{n}^{*}, f_{n}^{*})\}$, and an element   $(u^*,f^*)\in V\times L^2(D)$ such that (\ref{se1}--(\ref{se3}) hold, which concludes the proof. \hfill$\Box$

\begin{Corollary}	\label{cor4} Assume
	$(\ref{b})$--$(\ref{sm})$,  $(\ref{L3})$, $(\ref{L4})$ and $(\ref{L1m})$, $(\ref{L2m})$,    for all $\mu\ge 0$. Then, for each $\mu\ge 0$ there exists a pair
	$(u^*_\mu,f^*_\mu)$ such that
	\begin{equation}\label{3dp}
	(u^*_\mu,f^*_\mu)\in{\cal V}_{ad},\qquad \cL_\mu(u^*, f^*)=\min_{(u,f)\in {\cal V}_{ad}}\cL_\mu(u,f).
	\end{equation}
	Moreover, if $\{(u^*_n,f^*_n)\}$ represents a sequence pairs such that  $(u^*_n,f^*_n)$ is a solution of Problem $\cQ$ corresponding to $\mu_n\ge 0$ for each $n\in\mathbb{N}$, then
	there exists a subsequence of the sequence $\{(u_{n}^{*}, f_{n}^{*})\}$, again denoted by $\{(u_{n}^{*}, f_{n}^{*})\}$, and an element   $(u^*,f^*)\in V\times L^2(D)$ such that $(\ref{se1})$--$(\ref{se3})$ hold.	
\end{Corollary}

The proof of Corollary \ref{cor4} follows from arguments similar to those used in the proof of Corollary \ref{cor3} and, therefore, we skip it.

\begin{Corollary}	\label{cor5} Assume
$(\ref{b})$--$(\ref{sm})$, $(\ref{wK})$, $(\ref{G})$, $(\ref{L3})$, $(\ref{L4})$ and $(\ref{L1m})$, $(\ref{L2m})$,    for all $\mu\ge 0$. Moreover, assume that Problem $\cQ$ has a unique solution. Then Problem $\cQ$  is weakly  well-posed with respect to the Tykhonov triple  $(\ref{Ic})$--$(\ref{Cc})$.	
\end{Corollary}

\noindent	
{\it Proof.} Let $(u^*,f^*)$ be the unique solution to Problem $\cQ$ and let $\{u_n^*, f_{n}^*\}$ be
an approximating sequence for the Problem $\cQ$ with respect to the Tykhonov triple $(\ref{Ic})$--$(\ref{Cc})$ First, it follows from the proof of Theorem \ref{t2} that the sequence  $\{f_{n}^*\}$ is bounded in $L^2(D)$. Therefore, using arguments similar to those used in step i) of the proof of Theorem \ref{t1}   we deduce that the sequence $\{u_{n}^*\}$ is bounded in $V$. We conclude from here that the sequence $\{(u_n^*, f_{n}^*)\}$ is bounded in the product space $V\times L^2(D)$.
Second, a careful analysis of the proof of Theorem \ref{t2} reveals that $(u^*,f^*)$ is the weak limit  (in $V\times L^2(D)$) of any weakly convergent subsequence of the sequence   $\{(u_n^*, f_{n}^*)\}$. The two properties above allow us to  use a standard argument in order to deduce   that
the whole sequence
$\{(u_n^*, f_{n}^*)\}$ converges weakly in $V\times L^2(D)$ to  $(u^*,f^*)$, as
$n\to \infty$. Corollary \ref{cor5} is now a direct consequence of Definition \ref{def2} b).\hfill$\Box$

\medskip
We end this section with the following remarks.
First, Corollary \ref{cor3} provides a topological property of the set of solutions of Problem $\cQ$. Moreover, Corollary \ref{cor4} provides a continuity result of the solutions of this problem with respect to the parameter $\mu$ and, implicitely, with respect to the cost functional. Finally, Corollary \ref{cor5} provides a weakly well-posedness result for Problem $\cQ$, in the particular case when this problem has a unique solution. Such situations arise for specific boundary condition and cost functionals, as explained in \cite{BT1,BT2}.

\section{Some relevant examples }\label{s5}
\setcounter{equation}0

We start this section with examples of perturbed heat transfer problems which lead to hemivariational inequalities
of the form (\ref{hvp}) for which conditions $(\ref{wK})$  and $(\ref{G})$ are satisfied.

\medskip\noindent
{\bf A unilateral problem with penalty conditions on $D$.} The first problem we introduce in this section is obtained by considering a penalty version of the unilateral condition in (\ref{d1}). The problem is formulated as follows.

\medskip\noindent{\bf Problem $\cH_0$}. {\it
	Find a temperature field $u:D\to\R$ such that}
\begin{eqnarray}
&&\label{d10}\Delta u+f=\frac{1}{\lambda}\,p_0(u)\qquad{\rm a.e.\ in\ }D,\\ [2mm]
&&\label{d20}u=0\hspace{30mm}{\rm a.e.\ on\ }\Gamma_1,\\ [2mm]
&&\label{d30}u=b\hspace{30mm}{\rm a.e.\ on\ }\Gamma_2,\\ [2mm]
&&\label{d01}-\frac{\partial u}{\partial\nu}\in\partial j_\nu(u)\hspace{14mm}{\rm a.e.\ on\ }\Gamma_3.
\end{eqnarray}

Here $\lambda>0$ is a penalty parameter and $p_0$ is a function which is assumed to satisfy the following condition.

\begin{equation}
\left\{ \begin{array}{ll}
p_0\colon D\times\mathbb{R}\to\mathbb{R}
\ \mbox{is such that} \\ [1mm]
{\rm (a)}\  |p_0(\bx,r)-p_0(\bx,s)|\le L_0 |r-s|\\
\qquad \quad \mbox{for all} \,
r,\,s\in \mathbb{R},\ {\rm a.e.\ } \bx\in D,\ {\rm with}\ L_0>0,\\ [1mm]
{\rm (b) \ }
(p_0(\bx,r)-p_0(\bx,s))\,(r-s)\ge 0\\
\qquad\quad \mbox{for all}\ \
r,\,s\in \mathbb{R},\ {\rm a.e.\ } \bx\in D,\\ [1mm]
{\rm (c)} \ \bx\mapsto p_0(\bx, r) {\rm\ is\ measurable\ on}\ D \
{\rm for\ all \ }r\in \mathbb{R},\\[1mm]
{\rm (d)}\ p_0(\bx,r)=0 \ \ \mbox{if and only if} \ \ r \ge 0,
\ {\rm a.e.}\ \bx\in D.
\end{array}\right.
\label{p0}
\end{equation}

\noindent
A typical example of such function is given by
\begin{equation*}
p_0(x,r)=-cr^-\qquad {\rm for\ all \ }r\in \mathbb{R},\ \bx\in D
\end{equation*}
where $c>0$ and $r^-$ represents the negative part of $r$, i.e., $r^-=max\,\{-r,0\}$. Note that in this case the term $p_0(u)$ in (\ref{d10})  vanishes in the points of $D$ where $u\ge 0$ and equals
$c\,\frac{u}{\lambda}$ in the points of $D$ where $u<0$.

Using standard arguments it is easy to see that  Problem  $\cH_0$ leads to a variational formulation
of the form $(\ref{hvp})$
in which
\begin{eqnarray}
&&\label{wK0} \wK=\ \{\,v\in V\ :\ v=b\ \ {\rm on}\ \ \Gamma_2\,\},\\
&&\label{G0}\langle Gu,v\rangle=\int_{D}p_0(u)v\,da\qquad\forall\,u,\ v\in V.
\end{eqnarray}

We have the following result.

\begin{Proposition}\label{pr80}
	Assume $(\ref{b})$ and $(\ref{p0})$. Then, the set $(\ref{wK0})$ and the operator $(\ref{G0})$ satisfy conditions $(\ref{wK})$  and $(\ref{G})$, respectively.
\end{Proposition}

\noindent{\it Proof.}
 Conditions $(\ref{wK})$ and $(\ref{G})$(a) are clearly satisfied. Assume now that $u\in\wK$ and $v\in K$ where, recall, the set $K$ is defined by (\ref{K}).
Then, using (\ref{p0}) we see that
\begin{equation}\label{z53}
p_0(u)v\le 0\quad{\rm and}\quad p_0(u)u\ge 0\quad{\rm  a.e.\ in}\ D
\end{equation}
which imply that\ \
$p_0(u)(v-u)\le 0$\  \ a.e.\ in $D$. We conclude from here that
\[\int_Dp_0(u)(v-u)\,dx\le 0\]
and, therefore, condition $(\ref{G})$(b) holds.

Next, we assume that $u\in\wK$ and  $\langle Gu,v-u\rangle=0$ for all $v\in K$, which implies that
\begin{equation}\label{z33n}
\int_Dp_0(u)u\,dx=\int_Dp_0(u)v\,dx\qquad\forall\, v\in K.
\end{equation}
We now use inequalities (\ref{z53}) to deduce that
\begin{equation}\label{z33z}
\int_Dp_0(u)u\,dx=0.
\end{equation}
Therefore, (\ref{z53}), (\ref{z33n}) combined with  the implication
\begin{equation}\label{z63}
h\ge 0,\quad \int_{\omega} h=0\quad\Longrightarrow  \quad h=0\quad{\rm  a.e.\ \ on}\ \ \omega
\end{equation}
show that $p_0(u)u=0$\ \ ${\rm a.e.\ on\ } D$. This equality together with condition (\ref{p0})(d) implies that $u\ge 0$ ${\rm a.e.\ on\ } D$. Recall now that $u\in \wK$. We deduce from here that $u=b$ ${\rm a.e.\ on\ } \Gamma_2$. Therefore, $u\in K$, which concludes the proof.
\hfill$\Box$

\medskip\noindent
{\bf A unilateral problem with penalty conditions on $\Gamma_2$.} The second problem we introduce in this section is obtained by considering a penalty version of the boundary condition (\ref{d3}) on $\Gamma_2$. The problem is formulated as follows.

\medskip\noindent{\bf Problem $\cH_2$}. {\it
	Find a temperature field $u:D\to\R$ such that}
\begin{eqnarray}
&&\label{d11}u\ge 0,\qquad -\Delta u-f\ge 0,\qquad u(\Delta u+f)=0\qquad{\rm a.e.\ in\ }D,\\ [2mm]
&&\label{d21}u=0\hspace{36mm}{\rm a.e.\ on\ }\Gamma_1,\\ [2mm]
&&\label{d31}-\frac{\partial u}{\partial\nu}=\frac{1}{\lambda}\,p_2(u-b)\hspace{11mm}{\rm a.e.\ on\ }\Gamma_2,\\ [2mm]
&&\label{d41}-\frac{\partial u}{\partial\nu}\in\partial j_\nu(u)\hspace{19.5mm}{\rm a.e.\ on\ }\Gamma_3.
\end{eqnarray}

Here, again, $\lambda>0$ is a penalty parameter and  $p_2$ is a function which is assumed to satisfy the following condition.
\begin{equation}
\left\{ \begin{array}{ll}
p_2\colon \Gamma_2\times\mathbb{R}\to\mathbb{R}
\ \mbox{is such that} \\ [1mm]
{\rm (a)}\
|p_2(\bx,r)-p_2(\bx,s)|\le L_2 |r-s|\\
\qquad \quad \mbox{for all} \,
r,\,s\in \mathbb{R},\ {\rm a.e.\ } \bx\in\Gamma_2,
\ {\rm with}\ L_2>0, \\ [1mm]
{\rm (b) \ }
(p_2(\bx,r)-p_2(\bx,s))\,(r-s)\ge 0\\
\qquad\quad \mbox{for all}\ \
r,\,s\in \mathbb{R},\ {\rm a.e.\ } \bx\in\Gamma_2,\\ [1mm]
{\rm (c)} \ \bx\mapsto p_2(\bx, r) {\rm\ is\ measurable\ on}\ \Gamma_2 \
{\rm for\ all \ }r\in \mathbb{R},\\[1mm]
{\rm (d)}\ p_2(\bx,r)=0 \ \ \mbox{if and only if} \ \ r =0,
\ {\rm a.e.}\ \bx\in\Gamma_2.
\end{array}\right.
\label{p2}
\end{equation}

\noindent
A typical example of such function is given by
\begin{equation*}\label{p2n}
p_2(\bx,r)=cr\qquad {\rm for\ all \ }r\in \mathbb{R},\ \bx\in \Gamma_2
\end{equation*}
where $c>0$. 
Using standard arguments it is easy to see that  Problem  $\cH_2$ leads to a variational formulation
of the form $(\ref{hvp})$
in which
\begin{eqnarray}
&&\label{wK2} \wK=\ \{\,v\in V\ :\ v\ge 0\ \ {\rm in}\ \ D\,\},\\
&&\label{G2}\langle Gu,v\rangle=\int_{\Gamma_2}p_2(u-b)v\,da\qquad\forall\,u,\ v\in V.
\end{eqnarray}

We have the following result.

\begin{Proposition}\label{pr82}
Assume $(\ref{b})$ and $(\ref{p2})$. Then, the set $(\ref{wK2})$ and the operator $(\ref{G2})$ satisfy conditions $(\ref{wK})$  and $(\ref{G})$, respectively.	
\end{Proposition}

\noindent{\it Proof.}
Conditions $(\ref{wK})$ and $(\ref{G})$(a) are clearly satisfied. Assume now that $u\in\wK$ and $v\in K$.
Then, using (\ref{p2}) we see that\  \ $p_2(u-b)(v-u)= p_2(u-b)(b-u)$\ a.e. on $\Gamma_2$ and, using the properties (\ref{p2}) of the function $p_2$ we deduce that
$p_2(u-b)(v-u)\le 0$\  \ a.e.\ on $\Gamma_3$. We conclude from here that
\[\int_{\Gamma_2}p_2(u-b)(v-u)\,da\le 0\]
and, therefore, condition $(\ref{G})$(b) holds.

Assume now that $u\in\wK$ and  $\langle Gu,v-u\rangle=0$ for all $v\in K$. This implies that
\begin{equation*}\label{z33p}
\int_{\Gamma_3}p_2(u-b)(b-u)\,da=0.
\end{equation*}
Now, since (\ref{p2}) guarantees that $p_2(u-b)(b-u)\le 0$ a.e. on $\Gamma_2$ we use implication  (\ref{z63}), again, to deduce that $p_2(u-b)(b-u)=0$ a.e. on $\Gamma_2$. In now follows from condition
(\ref{p2})(d) that $u=b$ a.e. on $\Gamma_3$ and, therefore, $u\in K$.
\hfill$\Box$

\medskip\noindent
{\bf  A fully penalty problem.} The third problem we introduce in this section is obtained by considering a penalty version of both the conditions  (\ref{d1}) and  (\ref{d3}). The problem is formulated as follows.

\medskip\noindent{\bf Problem $\cH_{02}$}. {\it
	Find a temperature field $u:D\to\R$ such that}
\begin{eqnarray}
&&\label{d1p} \Delta u+f=\frac{1}{\lambda}\,p_0(u)\qquad\ {\rm a.e.\ in\ }D,\\ [2mm]
&&\label{d2p}u=0\hspace{32mm}{\rm a.e.\ on\ }\Gamma_1,\\ [2mm]
&&\label{d3p}-\frac{\partial u}{\partial\nu}=\frac{1}{\lambda}\,p_2(u-b)\hspace{7mm}{\rm a.e.\ on\ }\Gamma_2,\\ [2mm]
&&\label{d4p}-\frac{\partial u}{\partial \nu}\in\partial j_\nu(u)\hspace{15mm}{\rm a.e.\ on\ }\Gamma_3.
\end{eqnarray}

Here $\lambda>0$ is a penalty parameter and $p_0$ and $p_2$ are  functions which satisfy  conditions (\ref{p0}) and (\ref{p2}), respectively.
Using standard arguments it is easy to see that  Problem  $\cH_2$ leads to a variational formulation
of the form $(\ref{hvp})$
in which
\begin{eqnarray}
&&\label{wK02} \wK=V,\\
&&\label{G02}\langle Gu,v\rangle=\int_{D}p_0(u)v\,da+\int_{\Gamma_2}p_2(u-b)v\,da\qquad\forall\,u,\ v\in V.
\end{eqnarray}

We have the following result.

\begin{Proposition}\label{pr802}
Assume $(\ref{b})$, $(\ref{p0})$ and $(\ref{p2})$. Then, the set $(\ref{wK02})$ and the operator $(\ref{G02})$ satisfy conditions $(\ref{wK})$  and $(\ref{G})$, respectively.
\end{Proposition}

\noindent{\it Proof.}
Conditions $(\ref{wK})$ and $(\ref{G})$(a) are clearly satisfied. Assume now that $u\in\wK$ and $v\in K$.
Then, using (\ref{p0}) and (\ref{p2}) it is easy to see that $p_0(u)(v-u)\le 0$\  \ a.e.\ in $D$
and $p_2(u-b)(v-u)\le 0$\  \ a.e.\ on $\Gamma_3$.
We conclude from here that
\[\int_Dp_0(u)(v-u)\,dx+\int_{\Gamma_2}p_2(u-b)(v-u)\,da\le 0\]
and, therefore, condition $(\ref{G})$(b) holds.

Assume now that $u\in\wK$ and  $\langle Gu,v-u\rangle=0$ for all $v\in K$. This implies that
\begin{equation}\label{z33x}
\int_Dp_0(u)v\,dx+\int_{\Gamma_3}p_2(u-b)(b-u)\,da=\int_Dp_0(u)u\,dx.
\end{equation}
Now, recall that  (\ref{p0}) and  (\ref{p2}) guarantee that
\begin{eqnarray}
&&\label{z33y}\mbox{$p_0(u)v\le 0$\ \  a.e. \ in\  $D$,\quad $p_2(u-b)(b-u)\le 0$\ \ a.e.\  on\ $\Gamma_2$},\\[2mm]
&&\label{z33t}\mbox{$p_0(u)u\ge 0$\ \ a.e.\ in\ $D$}.\end{eqnarray}
We now combine equality  (\ref{z33x}) with inequalities (\ref{z33y}) and (\ref{z33t}) to find that
\begin{equation}\label{aa}\int_Dp_0(u)u\,dx=0.
\end{equation}
Next, (\ref{z33t}), (\ref{aa}) and (\ref{z63}) imply that $p_0(u)u=0$ a.e. in $D$ and, using  condition (\ref{p0}) we find that
\begin{equation}\label{101}
u\ge 0\qquad{\rm a.e.\ in}\ D.
\end{equation}
We conclude from here that $p_0(u)=0$ a.e. on $D$ and, therefore, (\ref{z33x}) yields
\begin{equation}\label{102}
\int_{\Gamma_3}p_2(u-b)(b-u)\,da=0.
\end{equation}
 Next, (\ref{z33y}), (\ref{102}) and (\ref{z63}) imply that $p_2(u-b)(b-u)=0$ a.e. on $\Gamma_2$. In now follows from condition (\ref{p2})(d) that
 \begin{equation}\label{103}
 u=b\qquad{\rm a.e.\ on}\ \Gamma_3.
 \end{equation}
 We now use (\ref{101}) and (\ref{103}) de deduce that $u\in K$ which concludes the proof.
\hfill$\Box$

\medskip

\medskip\noindent
{\bf  An example of cost functional.}  A large number of cost functionals for which our results in Section \ref{s4} hold can be considered. Here we restrict ourselves to provide the following example.  Assume that
\begin{eqnarray}
&&\label{r1}a_0>0,\qquad a_2>0,\qquad \phi\in L^2(\Gamma_2),\\[2mm]
&&\label{r2}\omega=[0,+\infty)\to\R\ \mbox{\rm is a continuous function such that}\ \omega(0)=1.
\end{eqnarray}	
Consider the cost functional ${\cal L}:V\times L^2(D)\to\R$ defined by
\begin{equation}
{\cal L}(u,f)=a_0\int_D f^2\,dx+ a_2\int_{\Gamma_2}(u-\phi)^2\,da
\end{equation}
and, for each $\mu\ge0$, let ${\cal L}_\mu:V\times L^2(D)\to\R$ be defined by
\begin{equation}
{\cal L}_\mu(u,f)=a_0\int_D f^2\,dx+ a_2\int_{\Gamma_2}(u-\omega(\mu)\phi)^2\,da
\end{equation}
Then, it is easy to see that the following result holds.

\begin{Proposition}\label{pr9}
Under assumptions $(\ref{r1})$ and $(\ref{r2})$ the functionals ${\cal L}$ and ${\cal L}_\mu$ satisfy conditions $(\ref{L1})$, $(\ref{L2})$, $(\ref{L3})$, $(\ref{L4})$ and $(\ref{L1m})$, $(\ref{L2m})$, for each $\mu\ge0$.
\end{Proposition}

\medskip\noindent
{\bf Final remarks.} We end this section with the remark that	
 Propositions \ref{pr80}--\ref{pr802} allow us to apply Theorem \ref{t1} and  and Corollaries \ref{cor1}, \ref{cor2} in the study of the	corresponding boundary  value problems. In this way we obtain the unique solvability  of Problems $\cH_0$, $\cH_2$ and $\cH_{02}$, the weak-strong continuous dependence of their weak solutions with respect to $f\in L^2(D)$ and the convergence of these solutions to the weak solution of Problem $\cH$, as $\lambda\to 0$. Moreover, Proposition \ref{pr9} allows us to apply Theorem \ref{t2} and  Corollaries \ref{cor3}--\ref{cor5} in the study of the corresponding optimal control problems. This allows to obtain the solvability of these optimùal control problems, the sequential compactness of the sets of their solutions and their continuous dependence with respect to the parameter $\mu$.


\vspace{8mm}

\section*{Acknowledgement}

\indent This project has received funding from the European Union's Horizon 2020
Research and Innovation Programme under the Marie Sklodowska-Curie
Grant Agreement No 823731 CONMECH.

\vskip 6mm
\noindent
{\it Authors'address:}

\vskip 2mm
\noindent
{\it  Mircea Sofonea},  Laboratoire de Math\'ematiques et Physique
	University of Perpignan Via Domitia, 52 Avenue Paul Alduy, 66860 Perpignan, France, e-mail: sofonea@univ-perp.fr

\vskip 2mm
\noindent	
{\it  Domingo A. Tarzia},	Departamento de Matematica-CONICET, Universidad Austral,
		Paraguay 1950, S2000FZF Rosario, Argentina, e-mail: DTarzia@austral.edu.ar

\end{document}